\documentclass[11pt]{article}

\usepackage{amsmath,amssymb,amsthm,color,graphicx,url}
\usepackage{times}
\usepackage{euscript}

\usepackage{latexsym}

\def\appendix{\bigskip\par
        \def\section{
        \secdef\Appendix\sAppendix}
        \setcounter{section}{0}
        \setcounter{subsection}{0}
        \gdef\thesection{\Alph{section}}
        }

\newtheorem{theorem}{Theorem}[section]
\newtheorem{lemma}[theorem]{Lemma}

\begin{document}

\begin{titlepage}

\title{Towards an optimal algorithm for recognizing Laman graphs\thanks{The authors would like
to thank Ileana Streinu for organizing the 2005 and 2006 Barbados workshops,
partially supported by NSF. This work would not have been possible without
participation in those workshops.}}

\author{Ovidiu Daescu\thanks{
	Department of Computer Science,
	University of Texas at Dallas,
	Richardson, TX 75080, USA.
	E-mail: {\tt daescu@utdallas.edu}. Daescu's research is supported
in part by NSF award CCF-0635013.}
\and
	Anastasia Kurdia\thanks{
	Department of Computer Science,
	University of Texas at Dallas,
	Richardson, TX 75080, USA.
        E-mail: {\tt akurdia@utdallas.edu}}}

\date{}
\maketitle

\begin{abstract}

Laman graphs are fundamental to rigidity theory. In computational geometry, they
are closely related to {\em pointed pseudo-triangulations} of planar point sets through
a property that states that the underlying graphs of pointed pseudo-triangulations
are Laman graphs.
A graph $G$ with $n$ vertices and $m$ edges is a {\em Laman graph}, or equivalently
a generic minimally rigid graph, if $m=2n-3$ and every induced subset of $k$ vertices
spans at most $2k-3$ edges.

We discuss the problem of recognizing Laman graphs.
Specifically, we consider the {\bf Verification problem:} {\em Given a graph G with $n$
vertices, decide if it is Laman}.

The previously best known algorithm for the verification problem takes $O(n^{3/2})$ time.
In this work we present an algorithm that takes $O(T_{st}(n)+n \log n)$ time, where
$T_{st}(n)$ is the best time to extract two edge disjoint spanning trees from $G$
or decide no such trees exist. So far, it is known that $T_{st}(n)$ is $O(n^{3/2})$.
Our algorithm exploits a known construction called red-black hierarchy (RBH), that is a
certificate for Laman graphs.
Previous algorithms construct the hierarchy in $O(n^2)$ time.
Our contribution is two-fold. First, we show how to verify if $G$ admits an RBH in
$O(n \log n)$ time and argue this is enough to conclude whether $G$ is Laman or not.
Second, we show that the RBH can be actually constructed in $O(n \log n)$ time using a two
steps procedure that is simple and easy to implement.

Finally, we point out some difficulties in using red-black hierarchies to compute a Henneberg
construction, which seem to imply super-quadratic time algorithms when used for embedding
a planar Laman graph as a pointed pseudo-triangulation.

\end{abstract}

\end{titlepage}

\section{Introduction}

Generic minimally rigid graphs in the plane, also known as Laman graphs, are fundamental to
rigidity theory~\cite{motionplan,haas}.
A graph $G$ with $n$ vertices and $m$ edges is a {\em Laman graph}
if $m=2n-3$ and every induced subset of $k$ vertices
spans at most $2k-3$ edges.
In computational geometry, they
are closely related to {\em pointed pseudo-triangulations} of planar point sets through
a property that states that the underlying graphs of pointed pseudo-triangulations
are Laman graphs~\cite{motionplan}.
Thus, pointed pseudo-triangluations inherit the properties of Laman graphs.
For example, related to the work in this paper, it follows that if we double any edge of
a given pointed pseudo-triangulation then its underlying graph can be decomposed in two
edge disjoint spanning trees.
Moreover, while not all Laman graphs can be embedded as pointed pseudo-triangulations,
every {\em planar} Laman graph can be embedded as a pointed pseudo-triangulation~\cite{haas}.

In this paper we consider recognizing Laman graphs.
Specifically, we address the {\bf Verification problem:} {\em Given a graph G with
$n$ vertices, decide if it is Laman}.

Most existing verification algorithms take quadratic time
in the number of input vertices to recognize Laman graphs~\cite{rigidcomp,ber}.
A very elegant and simple algorithm is the {\em pebble game} algorithm, first proposed by
Jacobs and Hendrickson~\cite{jacobshend}, and generalized later on by Streinu, Lee, and Theran in a
number of papers~\cite{hypergraph,leepebble,graphdecomp,mapdecomp}.
The pebble game algorithm solves the verification problem in $O(n^2)$ time.

Recski~\cite{r-1984} and Lovasz and Yemini~\cite{ly-1982} proved that a graph $G=(V,E)$ is Laman
if and only if, for each edge $e \in E$, the multigraph $G \cup \{e\}$ is the union of two edge
disjoint spanning trees.
In the remaining of this section we assume an edge of $G$ has been doubled and $G$ denotes the
resulting graph.

A known subquadratic time algorithm is due to Gabow and Westermann~\cite{gabow} and requires
$O(n^{3/2})$ time.
They solve this problem in two steps:
(1) Find a 2-forest of $G$ (two edge disjoint spanning trees), which
is done in $O(n^{3/2})$ time, and
(2) Test if the top clump is empty: this is done in $O(n \log n)$ time and uses some
structures computed in step (1).
Thus, step (2) is coupled with step (1), in the sense that
if two edge disjoint spanning trees are given to step (2), computed by some arbitrary method,
then step (2) should be changed and could require asymptotically larger time.
Very recently, it was suggested to us that a method presented in~\cite{graphdecomp} can be
adapted to speed up the top clump test to $O(n)$ time, assuming the data structures computed
in step (1) are available.

A different verification algorithm was proposed recently by Bereg~\cite{ber}.
The method in~\cite{ber} performs a step-by-step decomposition of $G$, aiming to construct a
hierarchical decomposition $H$ of $G$, called a {\em red-black hierarchy} (RBH).
It is argued in~\cite{ber} that $G$ is a Laman graph if and only if it admits a RBH.
The RBH construction in~\cite{ber} has three steps: (1) Find two edge disjoint spanning trees,
by some method
(Bereg uses an $O(n^2)$ time algorithm to obtain the trees, but he could have used
the algorithm in ~\cite{gabow}, for $O(n^{3/2})$ time);
(2) Construct a red-black hierarchy, which is done in $O(n^2)$ time, and
(3) Certify the hierarchy, which is done in $O(n)$ time.
Since steps (2) and (3) do not depend on how step (1) is performed, Bereg's
method decouples the computation of the two edge disjoint spanning trees in step (1)
from the rest of the computation.
Let $T_{st}(n)$ be the time to find two edge disjoint spanning trees.
Step (1) takes $O(T_{st}(n))$ time,
step (2) takes $O(n^2)$ time~\cite{ber}, and step (3) takes $O(n)$ time, totaling
$O(T_{st}(n)+n^2)$ time.

We present an $O(T_{st}(n)+n \log n)$ time verification algorithm
based on the following simple observation: from Corollary 4 in~\cite{ber}, it is not necessary
to actually construct $H$ to decide $G$ is Laman; we only need to decide whether a RBH
decomposition $H$ {\em exists} for $G$.
Thus, steps (2) and (3) above from Bereg's algorithm become: (2) use the two spanning trees
to decide whether $G$ {\em admits} a RBH decomposition.

Our algorithm has two steps:
(1) Compute two edge disjoint spanning trees by the best possible
method.
We use the algorithm in~\cite{gabow} since this is the best we know
(if, say, a simple $O(n \log n)$ time algorithm is discovered for this part, we
will use that one). This step takes $O(n^{3/2})$ time.
(2) Given two edge disjoint spanning trees for $G$, we give a simple solution for deciding
whether $G$ admits a RBH decomposition, that uses depth-first search and segment trees only,
and takes $O(n \log n)$ time. This step is independent of how step (1)
is done. At the end of step (2) we know if $G$ is Laman or not.
Moreover, we also show that the RBH can be actually constructed in $O(n \log n)$ time
using a two steps procedure that is simple and easy to implement.
Thus, our algorithm decouples step (1) from step (2), achieving the
desirable feature of Bereg's method (to take advantage of future
improvements on step (1)), and solves the second step of the verification in $O(n \log n)$
time instead of $O(n^2)$ time.

Finally, we point out some difficulties in using red-black hierarchies to compute a Henneberg
construction, which seem to imply super-quadratic time algorithms when red-black hierarchies
are used for embedding a planar Laman graph as a pointed pseudo-triangulation.


\section{Red-black hierarchies}
\label{rbh}

Red-black hierarchies (RBH) are introduced in~\cite{ber} as follows.

A hierarchy $H(G,T_h,\alpha,\beta)$ for a given graph $G(E,V)$, $|V|=n$, is a
graph $H(E_h,V_h)$, $E_h=T_h \cup \beta (E)$. $T_h$ is a set of edges forming a rooted tree.
The function $\alpha: V \rightarrow L(T_h)$, defines a one-to-one correspondence between
the vertices of $V$ and the leaves of the tree, denoted as $L(T_h)$.
The function $\beta: E \rightarrow V(T_h)\times V(T_h)$ maps an edge $(u,v)$ of $G$
to the edge $\beta(u,v)=(\beta_1(u,v),\beta_2(u,v))$ of $H$ (called cross edge), so
that $\beta_1(u,v)$ and $\beta_2(u,v)$ are ancestors, but not common ancestors,
of $\alpha(u)$ and $\alpha(v)$, respectively.

A RBH is a hierarchy $H(G,T_h,\alpha,\beta )$ satisfying the following conditions:
(1) The root of the tree $T_h$ has exactly two children (root rule);
(2) A vertex is the only child of its parent if and only if it is a leaf  (leaf rule);
(3) For any cross edge its endpoints have the same grandparent but different parents in the
tree (cross-edge rule);
(4) Cross edges connect all grandchildren of a vertex and form a tree (tree rule).

\begin{figure}
\centerline{\includegraphics[width=3.5in]{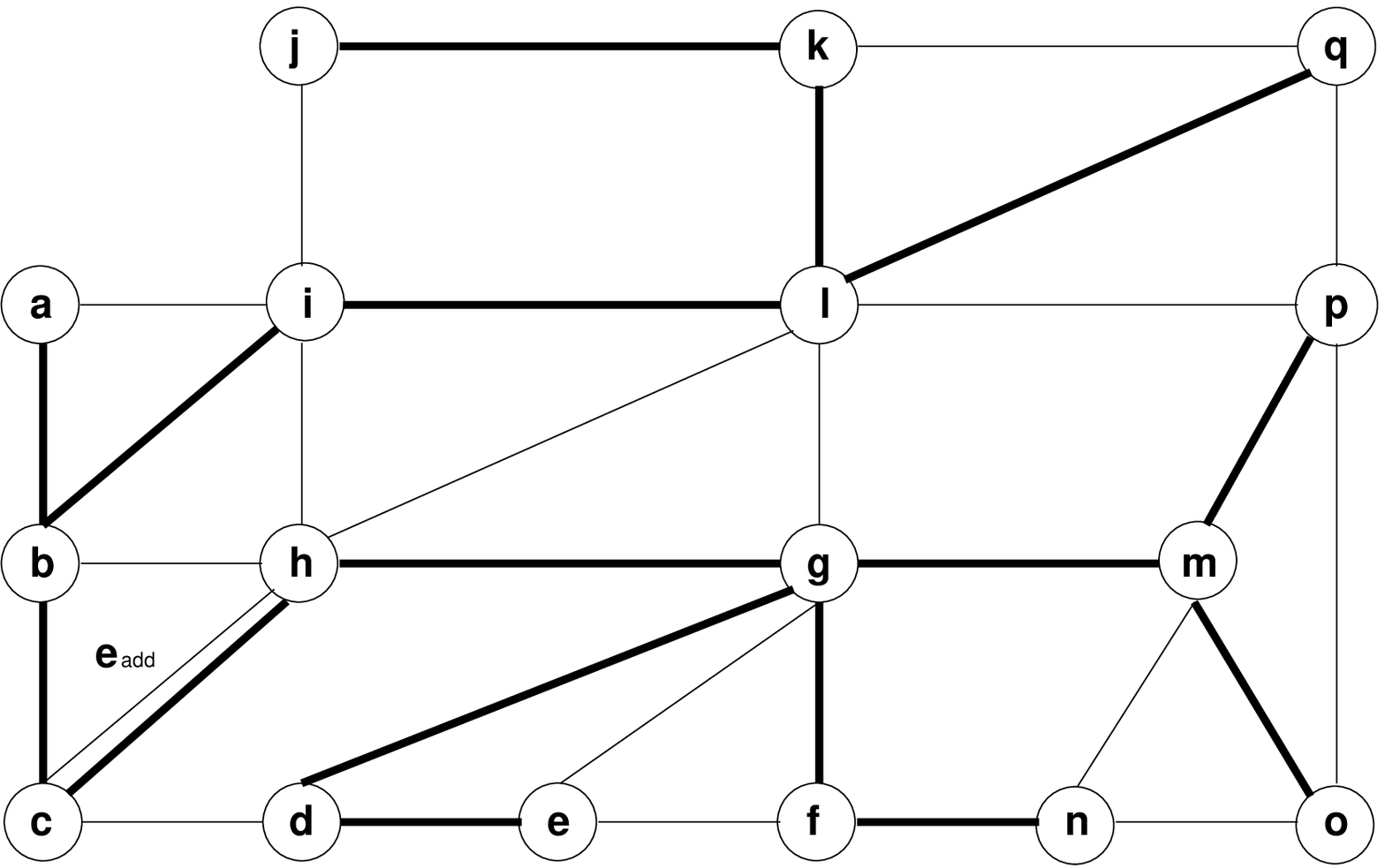}}
\caption{A graph $G^*=G \cup \{e_{add}\}$.}
\label{graph1}
\end{figure}

Given $G$, the construction of the RBH in~\cite{ber} has two major phases. First, a copy of an
edge of $G$, $e_{add}$, is added to $G$ and two edge-disjoint spanning trees,
$T^{r}$ (called {\em red tree}) and $T^{b}$ (called {\em black tree}),
are computed for $G^{*}=G\cup \{e_{add}\}$ using a known method (if no such trees exist,
then $G$ is not Laman and we stop). A graph $G^*$ and its two edge disjoint spanning trees
are shown in Figure~\ref{graph1}.
Second, a decomposition of $G^{*}$
is performed and a characterizing hierarchy $H=H(G^{*})$ is constructed in correspondence
with the steps of the decomposition.
We describe this decomposition~\cite{ber} below.

Suppose $e_{add} \in T^b$ and let $E(G^{*})=T^r \cup T^b$.
In the first step, a root $r_h$, corresponding
to $T^r$, is created in $H$ and is colored red. In the second step, the edge $e_{add}$
is removed from $T^b$ and two nodes corresponding to the resulting black trees $T^b_0$ and
$T^b_1$ are made children of $r_h$ in $H$ and are colored black.

Then, an iterative procedure is performed to construct $H$. At the end of step $i-1$, edges of
one color $c$ form a spanning forest $F^c=\{T^c_{0},\ldots, T^c_{l}\}$ of
$G=\{C_0,\ldots, C_l\}$,
where $C_i$ are connected subgraphs of $G$, and each element $T^c_{i}$ of $F^c$ is a
spanning tree of
its connected subgraph $C_i$. When restricted to these subgraphs, edges of the other color
$\overline{c}$ form a set
$F^{\overline{c}}=
\{
\{T^{\overline{c}}_{0,0}, \ldots, T^{\overline{c}}_{0,k_0}\},
\ldots,
\{T^{\overline{c}}_{l,0}, \ldots , T^{\overline{c}}_{l,k_l}\}
\}$.
Each element
$F^{\overline{c}}_{i}$ of $F^{\overline{c}}$ is a forest  spanning its respective
connected subgraph $C_i$.
The trees $T^{\overline{c}}_{i,j}$ are "linked" together in $G$
only with edges of color $c$.
There are $l+1$ vertices of color $c$ at the last level of $H$, each corresponding to
a tree from $F^c$.
At the beginning of the $i$-th step, all edges of color $c$ crossing the multi-cut defined by
$F^{\overline{c}}$ are found
and deleted from $G$.
At this point, $G$ consists of $\sum_{i=0}^{l}|F_i^{\overline{c}}|$ connected
subgraphs, and the trees $T^{\overline{c}}_{i,j}$ of color $\overline{c}$ are the spanning
trees of their respective subgraphs.
For each vertex $v_{h}$ of $H$ corresponding
to a tree $T^c_{i}$ from $F^c$, $k_i+1$ vertices of color $\overline{c}$ corresponding
to the trees
from $F^{\overline{c}}_{i}$ are created in $H$ as children of $v_{h}$.
A cross edge is added between
two vertices $x_{h}$ and $y_{h}$ at this now lowest level if the corresponding trees
$T^{\overline{c}}_{x}$ and $T^{\overline{c}}_{y}$ were previously in the same connected
component and got
separated at step $i-1$ by removing an edge of color $\overline{c}$ between them.
Addition of cross edges to $H$ completes the $i$-th step of the decomposition.

At step $i+1$, these actions are repeated for the swapped colors.
When at some step $j$ a cut of color $c$ is to be found and some connected component $C$ does
not have such cut, a vertex $l_h$ corresponding to a tree of color $c$ that spans $C$ is
created in $H$ and the decomposition stops for $C$.
The decomposition of $G$ ends when it has ended for all connected components.

The resulting graph $H$ characterizing the decomposition of the graph in Figure~\ref{graph1}
is shown in Figure ~\ref{hierarchy1}.
The entire decomposition process is given in appendix due to lack of space.

\begin{figure}
\centerline{\includegraphics[width=3.5in]{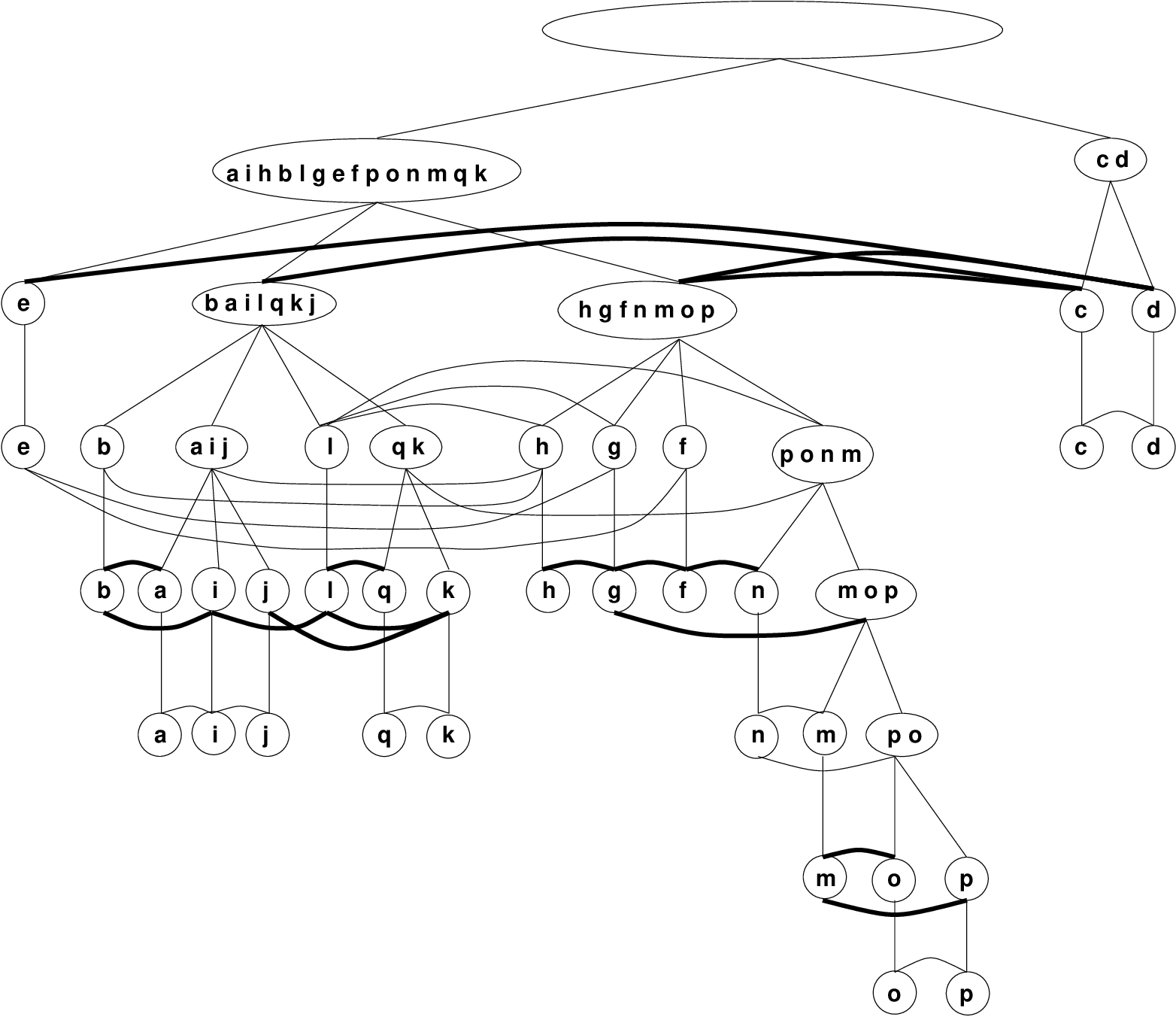}}
\caption{$H$ at the end of step $i=8$ that completes the decomposition of $G^*$.}
\label{hierarchy1}
\end{figure}

After $H$ is constructed, a check of whether $H$ satisfies the definition of the RBH is
performed. If the answer is positive, $H$ is a RBH and the corresponding graph $G$ is a
Laman graph. The method in~\cite{ber} takes $O(n^2)$ time to construct the
decomposition-characterizing graph $H$ and $O(n)$ time to verify that it satisfies the
definition of the RBH.

\section {A sufficient condition}
\label{sufficiency}

We show that if all edges are removed from $G$ during the decomposition process, the graph $H$
constructed from the decomposition is
always a RBH and thus $G$ is a Laman graph.

Vertices of $H$ correspond to spanning trees of connected subgraphs of $G$.
Anything marked by the
subscript $h$ in what follows refers to $H$.
$C(v_h)$ denotes the connected subgraph of vertex $v_h$.
$V(v_h)$ denotes the set of vertices of $C(v_h)$.
$T(v_h)$ denotes a tree spanning $V(v_h)$. $V(T)$ denotes
the set of vertices of $G$ spanned by the tree $T$.


We first prove that the four RBH rules introduced earlier always hold for the
decomposition-characterizing graph $H$ of
\textit{any} graph $G^{*}=G \cup e_{add}$, if the edge set of  $G^{*}$ can be partitioned into two
edge-disjoint spanning trees.
Let $color(v_h)$ denote the color associated with node $v_h$ (red or black).
If $c=color(v_h)$ is red then $\overline{c}$ is black and vice versa.

{\em Root rule}. At the very first step, $H$ is empty and a node $r_h$ of color $c$,
corresponding to the spanning tree $T^c$ that
does not contain the added edge $e_{add}$, is created in $H$.
The node $r_h$ is the root of $H$.
Then, $e_{add}$ is deleted from the other tree $T_{\overline{c}}$, which necessarily
creates exactly two trees of color $\overline{c}$ in $G^{*}$ and exactly two nodes
of color $\overline{c}$ in $H$ that are children of $r_h$, corresponding to these two trees.
Thus, the root rule always holds.

{\em Leaf rule}.
We first prove that if a vertex $v_h$ is the only child of its parent then $v_h$ is a leaf.
If a vertex $v_h$ is the only child of its parent, the connected subgraph $C(parent(v_h))$
could not be split any further during decomposition and $V(v_h)=V(parent(v_h))$.
At the step when $v_h$ was created, the decomposition process has stopped for $C(v_h)$:
there was just one tree of color $color(v_h)$ in $C(parent(v_h))$
and just one tree of color $\overline{color(v_h)}$
(otherwise $C(parent(v_h))$ would have been partitioned further and $v_h$ would have siblings).
Hence, the vertex $v_h$
corresponding to  $C(v_h)$ is a leaf in $H$.

Next, we prove that a leaf vertex cannot have any siblings.
Suppose there is a vertex $y_h$ having $k>1$ children and at
least one of them is a leaf. The vertex $y_h$ corresponds to a connected subgraph spanned
by a tree of
color $c=color(y_h)$ and a spanning forest of $k$ trees of color $\overline{c}$.
Each of its children
$x_h^i$ corresponds to a connected subgraph $C_i$, $i=1,2,\ldots,k$, spanned by a tree of color
$\overline{c}$ and a forest of color $c$ (possibly containing only one tree).
If this spanning forest
contains more than one tree, there are edges of color $\overline{c}$ in
$C_i$ connecting the trees of
the spanning forest. At the next step of the decomposition these edges
will be deleted, the spanning tree of color $\overline{c}$ will split into at least two
different trees and
corresponding vertices will be created in $H$ as children of $x_h^i$. Hence, $x_h^i$
cannot be a leaf vertex.
If the spanning forest of $C(x_h^i)$ contains just one tree then a vertex corresponding to that
tree, of color $c$, is created in $H$ as a child of $x_h^i$ and again, $x_h^i$
cannot be a leaf vertex.
This argument holds for every child of $y_h$, contradicting the assumption that at least one child
of $y_h$ is a leaf.

{\em Cross-edge rule}.
A cross edge is added between any two vertices $u_h$ and $v_h$ at step $i$ if their
corresponding vertex sets $V(u_h)$ and $V(v_h)$ previously belonged to one connected
subgraph $C_{u,v}$
and got separated at step $i-1$ by removing the edge between them. At level $i-2$ of $H$
there is always a vertex that corresponds to $C_{u,v}$.  The vertices at the same level of $H$
correspond to connected subgraphs that are disjoint subgraphs of $G$.
Hence, no other vertex at
level $i-2$ of $H$ can correspond to a connected subgraph containing $V(u_h)$, $V(v_h)$,
their subsets, or the union of their subsets. The vertex corresponding to the
connected subgraph $C_{u,v}$ is a common grandparent of  $u_h$ and $v_h$.

Again, according to the construction rules, parents of  $u_h$ and $v_h$ in $H$ correspond to
different connected subgraphs, so $u_h$ and $v_h$ have different parents.

{\bf Tree rule}.
If $k$ edges are removed from the tree $T$ spanning the vertex set $V(v_h)$ corresponding to
some vertex $v_h$ of $H$, $k+1$ new trees result from $T$ and $k+1$ nodes are created as
grandchildren of $v_h$ in $H$. For each edge $e$ deleted from $T$, a cross edge is added between
the vertices corresponding to the two sub-trees of $T$ that were connected by $e$. Each grandchild of
$v_h$ gets a cross edge incident to it. There are $k+1$ grandchildren of $v_h$ and $k$ edges
connecting them. The cross edges form a tree spanning all the grandchildren of $v_h$.

We have shown that red-black hierarchy rules hold for any $H$.
Then, we only need to check if $H$ satisfies the
general definition of a hierarchy.

\begin{lemma}
If all edges are removed from $G$ during the decomposition process then the characterizing graph
$H$ of $G$ satisfies the definition of hierarchy.
\end{lemma}

\noindent \textbf{Proof.}
The edges of $H$ are the union of the edges of the rooted tree $T_h$ and the cross edges.
There is a cross edge $e_h=(u_h,v_h)$ in $H$  for each edge $e=(u,v)$ of $G$:
The edge $e$ is deleted from $G$ when it crosses the cut separating $u$ from $v$; according
to the construction rules, a cross edge is then added between the vertices of $H$
corresponding to the connected components of $u$ and $v$ at the current step.

There is one-to-one correspondence between the leaves of $T_h$ and the vertices of $G$:
Since the graph splitting
procedure continues until all edges are removed, each vertex $v$ of $G$
is eventually disconnected from
the rest of the graph by deleting an edge of some color $c$. A vertex $l_h$ corresponding to a
tree of color $c$ spanning the connected subgraph $C_v=\{v\}$ is then created in $H$. Since $C_v$
cannot be split further, the decomposition stops for $C_v$ and $l_h$
becomes a leaf vertex of $T_h$.
Also,
there is no leaf vertex in $H$ that does not correspond to a vertex of $G$. Suppose there exists
such a vertex in $H$. Then, it corresponds to a tree spanning a
connected subgraph $C_{x}$, with $|C_x|>1$. This
means that $C_x$ contains edges that were not deleted during the decomposition of $G$,
a contradiction.

The endpoints of $e_h$ are ancestors of $\alpha(u)$ and $\alpha(v)$, respectively,
but they are not their common ancestors:
recall that $\alpha(u)$ and $\alpha(v)$ are the leaf vertices of $H$ corresponding
to vertices $u$ and $v$ of $G$. The vertex $u_h$ corresponds to some connected component $C(u_h)$.
The vertices of $H$ that are descendants of $u_h$ correspond to connected components over the
subsets of $V(u_h)$. The leaf vertices of $H$ that correspond to vertices in $V(u_h)$ are the
descendants of $u_h$ in $H$. Since $u \in V(u_h)$, $u_h$ is an ancestor of $\alpha(u)$.
Similarly, $v_h$ is an ancestor of $\alpha(v)$. The vertices  $u_h$ and $v_h$ are connected
by a cross edge corresponding to two disjoint connected subgraphs.
Therefore, $u_h$ cannot be an ancestor of $\alpha(v)$, and the end vertices of $e_h$ are not
common ancestors of $\alpha(u)$ and $\alpha(v)$.
\hfill $\Box$

\begin{lemma}{If $G$ has edges left at the end of the decomposition process, the
characterizing graph $H$ of $G$ does not satisfy the definition of hierarchy.}
\end{lemma}

\noindent \textbf{Proof.}
If there are non-deleted edges of $G$ when the decomposition stops, then there
are no corresponding edges for them in $H$. In addition, we do not have
a one-to-one map from $V$ to $L(T_h)$:
some leaves of $T_h$ correspond to connected subgraphs containing several vertices.
\hfill $\Box$

Thus, building $H$ is not required for certifying Laman graphs: just
decompose $G$ based on the rules in~\cite{ber} and check if
$G$ has edges left when the decomposition ends.

\section{The decomposition algorithm}
\label{groupingidea}

We have shown that building $H$ is not required for verifying that $G$ is a Laman graphs.
It is sufficient to perform the decomposition of $G$ according to the rules from~\cite{ber}
and then check whether there are edges left in $G$.
The decomposition algorithm has some notable features. At each step edges of only one color
are deleted. The groups of red and black edges are deleted in turns. At each step, except the
first and the last ones, at least one edge is deleted from $G$ (some edges may never be
deleted). Thus, the edges of $G$ can be grouped so that edges of one group are deleted from $G$ at
the same step.
The decomposition process provides a natural order on these groups.
We denote this ordering as $g=(g_2,g_{3},\ldots,g_{k})$, where the index $i$ of $g_i$
corresponds to the step at which the edges of the group $g_i$ are deleted.

Instead of $H$, we use $g$ to characterize the graph decomposition. Our main goal now is to
speed up the decomposition algorithm from~\cite{ber} using the following simple observation:
deletion of any edge $e=(u,v)$ from its tree (of color $color(e)$), where $u$ is a parent of
$v$ in a DFS ordering
of the tree of color $color(e)$, always forms two trees such that one of them is rooted at
$v$ and all nodes in that tree are descendants of $v$.

We slightly modify the graph decomposition algorithm from~\cite{ber}.
The edges to be deleted
at the next step are identified at the end of the preceding step and marked for deletion.
At the first step, $e_{add}$ is marked for deletion
(and no other action is performed).
Each iterative step in $G$ consists of removing the marked edges of some color $c$
and identifying and marking
the edges crossing the cuts of the opposite color $\overline{c}$ that appear after removing
the marked edges. We also note that once the original graph has split into several connected
subgraphs, the decomposition proceeds independently on each subgraph, and the problem of
finding the edges to be deleted at the subsequent step can be viewed as several independent
subproblems, each on a distinct connected subgraph.

Consider the graph $G^*$ and its two edge disjoint spanning trees
$T^c$ and $T^{\overline{c}}$,
rooted at vertices $r^c$ and $r^{\overline{c}}$, respectively.
Let $DFS(c)$ be the depth-first
search traversal of $G^*$ starting at $r^c$ and using only edges of color $c$, where $c$ is
either red or black.
We assign each vertex of $G$ two DFS order numbers, one from $DFS(red)$ and another one
from $DFS(black)$.
New edges are never added to the trees, so the numbers never change.
For any edge of color $c$, it is always possible
to establish the parent-child relationship of its endpoints by looking at their DFS numbers
for color $c$.
Whenever an edge $e$ is mentioned in this text as a vertex pair, the first vertex is always the
parent of the second vertex in $DFS(color(e))$.

When an edge $e=(u,v)$ of color $c$
is deleted from a
tree $T^c_{k}$ rooted at some $r^c$ and spanning a connected subgraph $C_k$, two trees emerge:
$T^c_{i}$ rooted at $r^c$ and $T^c_{j}$ rooted at $v$. Only the vertices of $T^c_{j}$ are
descendants of $v$ in DFS(c).
The ancestor/descendant relationship can be established in the $DFS(c)$
tree by looking at the discovery and finish times ($d^c[\cdot]$ and $f^c[\cdot]$, respectively)
of the vertices.

\begin{lemma}
\label{lem-cut}
{An edge $(x,y)$ of color $\overline{c}$ crosses the cut $(V(T^c_{i}), V(T^c_{j}))$
induced by the deletion of the edge $(u,v)$ of color $c$ if and only if
one of its endpoints is a descendant of $v$ and the other one is not, i.e., exactly one of its
endpoints discovery times is in $t=[d^c[v],f^c[v]]$.}
\end{lemma}

\noindent \textbf{Proof.}
If $d^c[x] \not\in t$ and $d^c[y] \in t$, then $x \in T^c_{i}$ and $y \in T^c_{j}$, so $e$
clearly crosses the cut. A symmetric argument applies if $d^c[x]  \in t$ and $d^c[y] \not\in t$.

If $d^c[x] \not\in t$ and  $d^c[y]\not\in t$, neither $x$ nor $y$ are in $T^c_{j}$, so both
endpoints of $e$ are in $T^c_{i}$ and $e$ does not cross the cut.
If $d^c[x] \in t$ and $d^c[y] \in t$, both endpoints of $e$ are in $T^c_{j}$ and $e$ does not
cross the cut.
\hfill $\Box$

From Lemma~\ref{lem-cut} it follows that if we associate an interval $[d^c[u],d^c[v]]$
with every edge of color $\overline{c}$, the intervals corresponding to the edges crossing
the cut have exactly one endpoint in $t$.

We identify such intervals using a segment tree data structure
enhanced with two
lists at each internal node, one sorted by the start time of the intervals stored at the node
and one sorted by their finish time.
A segment tree~\cite{prep} is a
balanced binary search tree that stores a set of intervals with endpoints from a finite set
of abscissae (intervals corresponding to edges of color $\overline{c}$, for example).
Each of its nodes $u$ has
an interval $I(u)$ associated with it and stores a list of input intervals intersecting $I(u)$.
Binary search in a segment tree allows to report the intervals containing a query point.

In our case, the endpoints of the intervals are integer numbers, so an interval containing a
point $p\pm\Delta$, for any $0<\Delta<1$, contains the point $p$ as well.
First, we find the intervals
with one endpoint before  $d^c[v]$ and the other endpoint in $t$ by querying for intervals
containing the point $d^c[v]-\Delta$. Second, we find the intervals with one endpoint in $t$
and the other endpoint after $f^c[v]$ by querying for intervals containing the
point $f^c[v]+\Delta$.

To ensure that each returned interval has an endpoint in $t$ we augment the standard segment
tree by storing two sorted lists at each node, instead of just one list.
With each node $u$, we store a list $L_{finish}(u)$
of intervals that intersect $I(u)$ that is sorted by the finish time of the intervals in
non-decreasing order; similarly, the list $L_{start}$ stores the same intervals sorted by their
starting time in non-increasing order. We give both queries above an additional parameter:
$f^c[v]$ for the first one and $d^c[v]$ for the second one. The first query only looks at the
lists $L_{finish}$ and reports the intervals that
have their right endpoint no greater than $f^c[v]$.
The second query only looks at the lists $L_{start}$
and reports the intervals that have their starting point no later than $d^c[v]$.
Thus, this data structure allows us to return intervals with exactly one endpoint in $t$.
Each query with an edge (interval $t$) takes
$O(\log n+k)$ time, where $k$ is the number of intervals (crossing edges) reported.
To avoid reporting an interval more
than once, the interval is deleted from the segment tree (including the sorted lists
associated with the nodes that store it) when it is returned by a query.
This can be easily done in $O(\log n)$ time.
Having two segment trees, one for the red intervals $[ d^{red}[u], d^{red}[v] ]$ of the black
edges and the another one for the black intervals $[ d^{black}[u], d^{black}[v]]$ of the red edges,
allows to efficiently identify edges of the cuts at each step
of the decomposition.

\begin{lemma}{The decomposition of $G$ can be done in $O(n \log n)$ time.}
\end{lemma}

\begin{theorem}{Given a graph $G$ with $n$ vertices and $m$ edges deciding whether $G$ is a
Laman graph or not can be done in $O(T_{st}(n)+n \log n)$ time, where $T_{st}(n)$ is the time to
extract two edge disjoint spanning trees from $G$ or decide no such trees exist.}
\end{theorem}

\noindent \textbf{Proof.}
We can check that $m=2n-3$ in $O(n)$ time.
Finding two edge disjoint spanning trees or deciding no such trees exist takes
$T_{st}(n)$ time. The best known algorithm so far for this task has
$T_{st}(n)=O(n^{3/2})$ time~\cite{gabow}.
The decomposition takes $O(n \log n)$ time:
$O(n \log n)$ for the segment trees,
$O(n\log n )$ to answer all
queries, and $O(n)$ to check if $G$ has any edges left
at the end of the decomposition.
\hfill $\Box$

\section{The reconstruction algorithm}
The order in which edges are deleted from  $G$ during the decomposition determines the structure
of the corresponding red-black hierarchy $H$, so given  $g$, one can unambiguously construct $H$ in
top-down fashion according to the rules from \cite{ber}. The vertices of $H$ correspond to
subtrees of
$T^c$ and $T^{\overline{c}}$,
and there is a vertex in $H$ for each distinct sub-tree (of $T^c$ or $T^{\overline{c}}$)
that appeared during the decomposition of $G$.
In the original approach, to construct the $i$-th level of $H$, one has to know the spanning sub-trees
at step $i-2$ of the decomposition and to figure out what trees appear after removal of edges at the
beginning of step $i$. It takes $O(n)$ time to find the emerging trees.

We consider
the decomposition process in reverse order (i.e. start from $n$ red and $n$ black disjoint trees
and add edges to them until two spanning trees are formed), and take advantage of the fact that it
is faster to union the disjoint sets into larger sets than to partition the trees into
disjoint sub-trees.
As a result, the proposed bottom-up construction method is faster and produces the same graph
$H$ as the top-down approach.

\begin{figure}
\centerline{\includegraphics[height=1.5in]{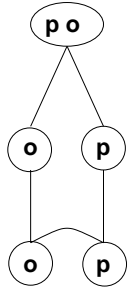}}
\caption{$H$ after considering edges of $g_8=\{(o,p)\}$.}
\label{reconstructionk}
\vspace{-0.1in}
\end{figure}

\begin{figure}
\centerline{\includegraphics[width=1in]{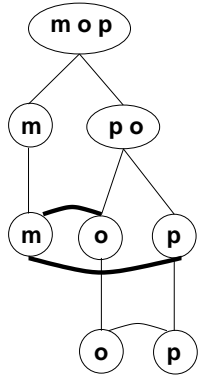}}
\caption{$H$ after considering edges of $g_8$ and $g_7=\{(m,o),(m,p)\}$.}
\label{reconstructionk-1}
\vspace{-0.1in}
\end{figure}

The last group $g_k$ of $g=(g_2,g_{3},\ldots,g_{k})$ contains edges of some color $c$ deleted at
the very last step of the decomposition. Each endpoint $v$ of edges of $g_k$ corresponds to a subtree
of $G$ of color $c$ spanning only the vertex $v$. A leaf node  $v_h=\alpha(v)$  is added to
the $k$-th
level of $H$ for each such vertex $v$. Only one leaf vertex is created for the endpoint shared by
multiple edges from $g_k$. For every edge $(u,v)$ of $g_k$ a corresponding cross edge
$\beta(u,v)=(\alpha(u), \alpha(v))$ is added to $H$.
For every leaf vertex $\alpha(v)$ of $H$, its parent should be at level $k-1$ of $H$,
corresponding to a subtree in $G$ that is of color $\overline{c}$ and spans only the vertex $v$.
Such parent vertex $v_h^p=parent(\alpha(v))$ is added to level $k-1$ of $H$ along with a
tree edge connecting $v_h^p$ and $\alpha(v)$ (we call this the {\em parent creation rule}).
The vertices of $H$
connected by a cross edge have the same grandparent. For every cross edge tree $T_k^j$ formed at
the $k$-th level, a vertex $v_{h}^g$ is added to level $k-2$ of $H$, as well as a tree edge
connecting $v_h^g$ and $v_h^p$, for every $v_h \in T_k^j$.
We have completed level $k$
of $H$ as well as added some elements to the two upper levels. See Figure~\ref{reconstructionk}
and Figure~\ref{reconstructionk-1} for an illustration.

At the $i$-th iterative step  for each cross edge $(x,y)$ of $g_{i}$ of color $c$
two vertices $v_h^x$ and $v_h^y$ on the $i$-th level of $H$ are identified.
They correspond to trees in $G$ of color $c$  that contained $x$ and $y$ respectively
at the $i$-th step of the decomposition.
If for some endpoint $x$ of an edge from $g_i$  $v_h^x$ does not exist on the $i$-th level of $H$,
a new vertex $v_h^x$ should be created at the $i$-th level and a parent for it should be added
following the parent creation rule.
Then the cross edge corresponding to $(x,y)$ is added to $H$ between $v_h^x$ and $v_h^y$.
After all edges of $g_{i}$ are considered, all cross-edges of the $i$-th level of $H$ are in place.
For each cross edge tree
$T_i^j$ formed at the $i$-th level of $H$, a node is added to level $i-2$ of $H$.
That grandparent node becomes a parent of the parents of the vertices of $H$ spanned by the
cross-edge tree $T_i^j$. At this time the $i$-th level of $H$ is complete and
levels $i-1$ and $i-2$ of $H$ are partially constructed.
Repeating these steps for all $g_i$, $i>2$, yields the RBH $H$.

\begin{lemma}{Given two edge-disjoint spanning trees for $G^*$, a red-black hierarchy for $G$,
if it exists, can be constructed in $O(n \log n)$ time.}
\end{lemma}

\noindent \textbf{Proof.}
Obtaining $g$ for $G^*$ takes $O(n \log n)$ time (Lemma 4.3). The time spent on reconstructing
one level is proportional to the number of cross edges at that level.
The total number of cross edges is $O(n)$. We use a standard UNION-FIND data structure for
maintaining the vertices of $G$ that the vertices of $H$ correspond to at each step
(notice that the actual trees of color $c$ or $\overline{c}$ defined by those vertices in $G$
are not needed to construct $H$). This allows to complete the reconstruction phase in
$O(n \log n)$ time, so the total time for constructing the RBH is $O(n \log n)$.

\hfill $\Box$

\vspace{-0.2in}

\section{RBH and Henneberg construction}

In this section
we point out some difficulties in using red-black hierarchies to compute a Henneberg
construction for $G$, which seem to imply super-quadratic time algorithms when red-black hierarchies
are used for embedding a planar Laman graph as a pointed pseudo-triangulation.

The vertex $e$ of the graph in Figure~\ref{graph1} falls in case 4 a) from \cite{ber}:
the grandparent of the leaf vertex $\alpha(e)$ has more than two children, $\alpha(e)$ has two
incident cross edges and its immediate parent has one incident cross edge
(see Figure~\ref{hierarchy1}). During the Henneberg construction we need to remove vertex $e$ from $G$ along with its incident
edges $(e,g)$, $(e,f)$, and $(e,d)$, and insert an edge between two of the vertices $g$, $f$, and
$d$ to restore the Laman property of the modified graph.
Note that $G$ already has edges $(g,d)$ and $(g,f)$, and thus adding the edge between $d$ and $f$
is the only option.

\begin{figure}
\centerline{\includegraphics[width=3.0in]{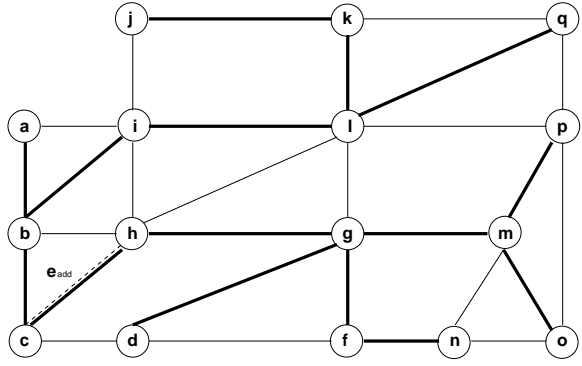}}
\caption{$G \setminus \{e\}$.}
\label{noegraph}
\vspace{-0.15in}
\end{figure}

\begin{figure}
\centerline{\includegraphics[width=3.5in]{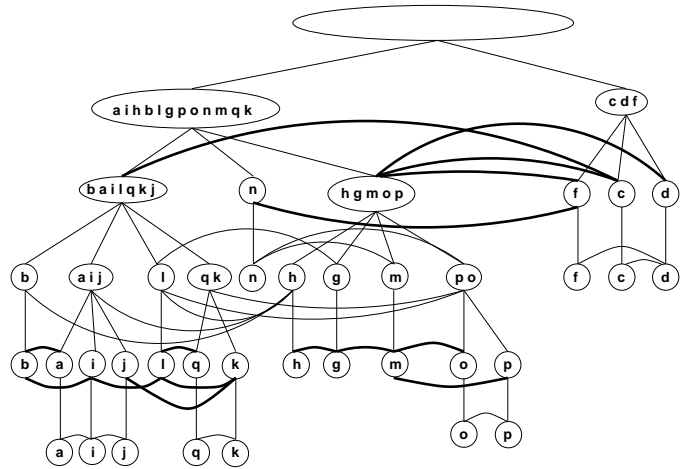}}
\caption{$H $ for $G \setminus \{e\}$.}
\label{noehierarchy}
\end{figure}

Having removed $e$ from $H$, we need to restore the properties of the hierarchy as if $e$ was
never present in the original graph. This does not appear to be an easy task, since the red-black
hierarchy for the graph $G \setminus \{e\}$ (Figure~\ref{noegraph}) differs significantly from the
one for $G$ (Figure~\ref{noehierarchy}).
Thus, it seems one would need to recompute the RBH for the resulting graph (starting with finding
two edge disjoint spanning trees), which would take $O(T_{st}(n)+n \log n)$ time rather than
$O(n)$ time in~\cite{ber}.
Essentially, obtaining the new RBH from the old one would be as difficult as obtaining two edge
disjoint spanning trees for the new graph from the edge disjoint spanning trees of the
original graph.

The argument above implies that,
over $O(n)$ steps, the Henneberg decomposition would take time $O(n(T_{st}(n)+n \log n))$.
Accordingly, embedding a planar Laman graph as a pointed pseudo-triangulation using red-black
hierarchies would require $O(n^{2.5})$ time using the best known algorithm for finding
two edge-disjoint spanning trees, which gives $T_{st}(n)=O(n^{1.5})$~\cite{gabow}.

\baselineskip=12pt
\small
\bibliography{bibliography}
\bibliographystyle{plain}

\newpage

{\bf\Large Appendix}

{\bf A sample construction of the decomposition characterizing graph}

\begin{figure}[h]
\centerline{\includegraphics[width=3.5in]{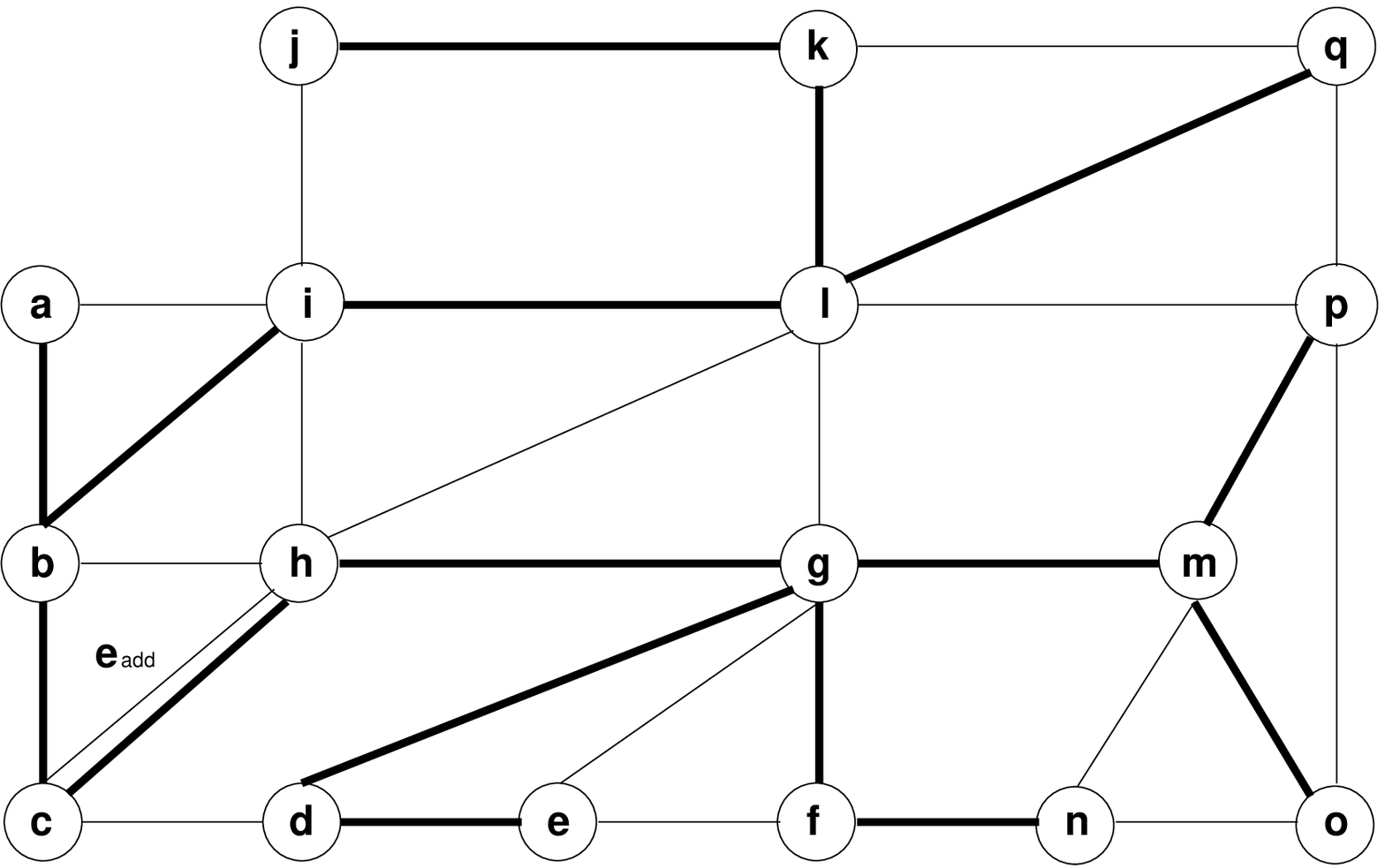}}
\caption{Original sample graph $G^*$.}
\label{graph}
\end{figure}

As an example we consider the decomposition of the graph from Figure~\ref{graph}.
The red tree (drawn with thick lines) is rooted at the vertex $b$ and the black tree is rooted at the vertex $a$.
All vertices of $H$, except the top one, are marked with the lists of vertices of their
corresponding subtrees. The edges of $G^{*}$ deleted at the $i^{th}$ step are shown with
dashed lines.

\begin{figure}
\centerline{\includegraphics[width=3.5in]{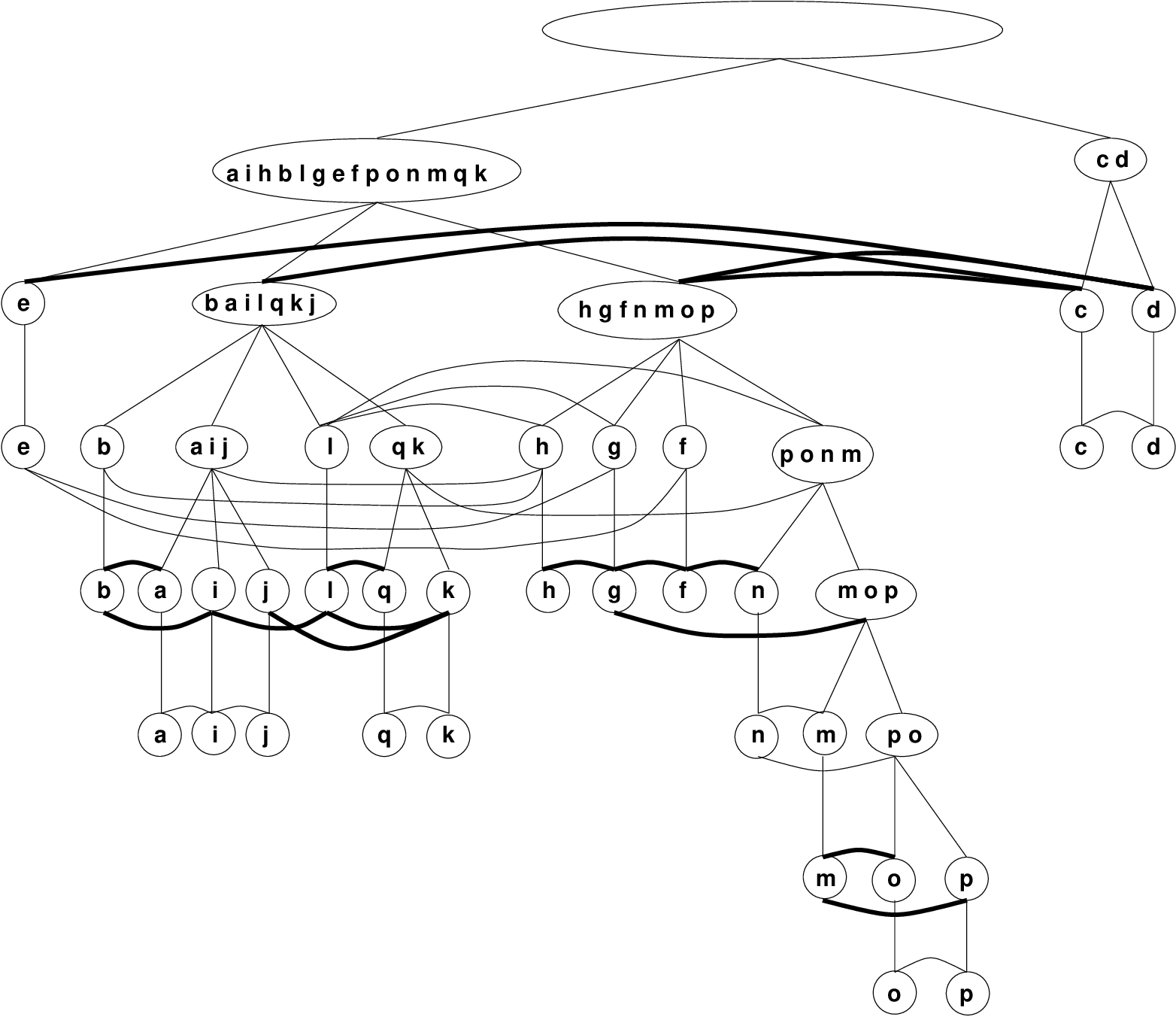}}
\caption{$H$ at the end of step $i=8$ that completes decomposition of $G^*$.}
\label{hierarchy}
\end{figure}

Figure~\ref{hierarchy} presents the resulting graph $H$ characterizing the decomposition of $G$.

\begin{figure}[h]
\centerline{\includegraphics[width=3.5in]{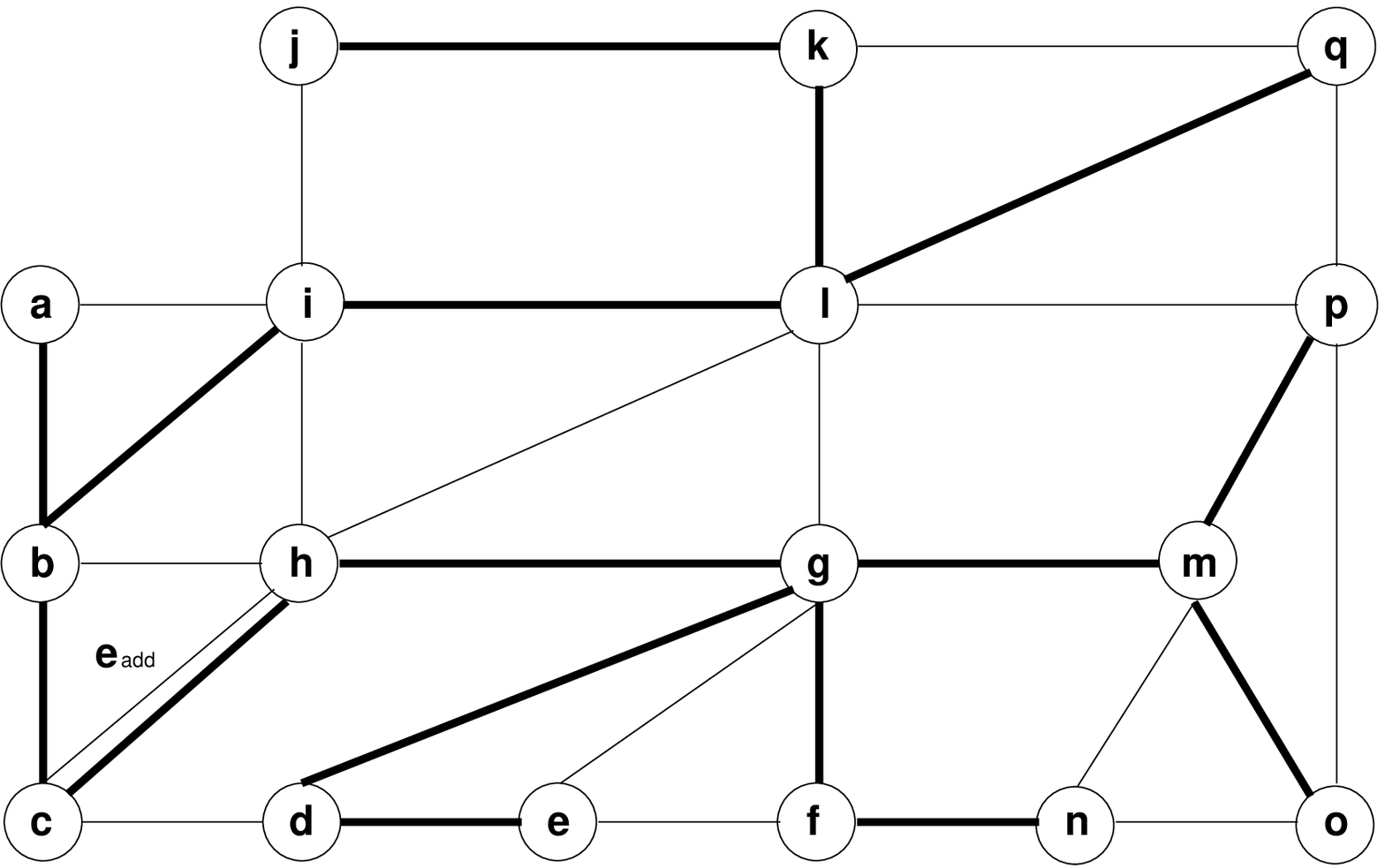}}
\caption{$G^*$ at the end of step $i=1$, no changes in the original graph.}
\end{figure}

\begin{figure}
\centerline{\includegraphics[width=3.5in]{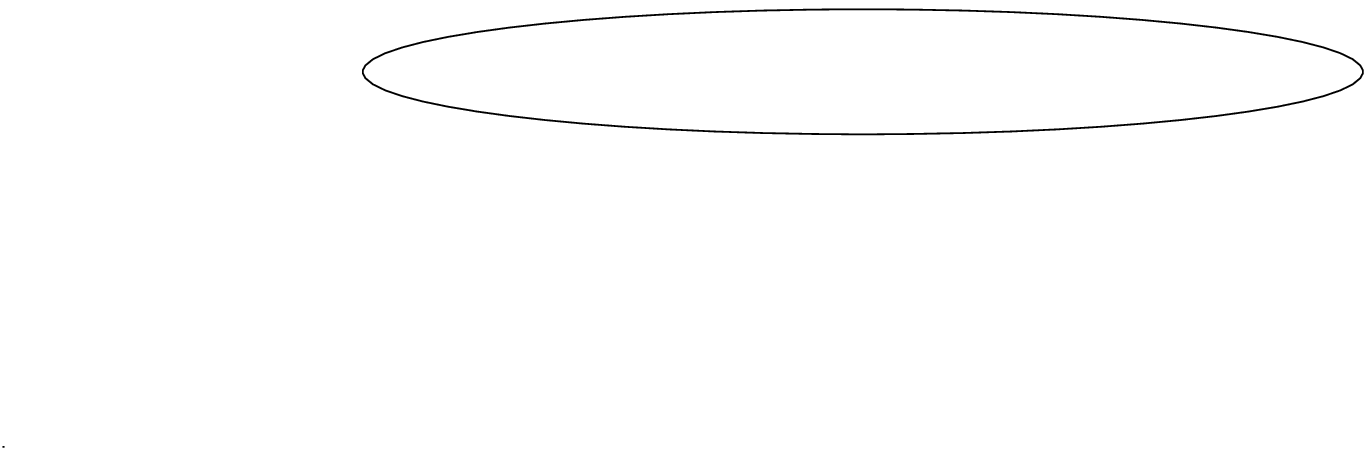}}
\caption{$H$ at the end of step $i=1$.}
\end{figure}

\begin{figure}
\centerline{\includegraphics[width=3.5in]{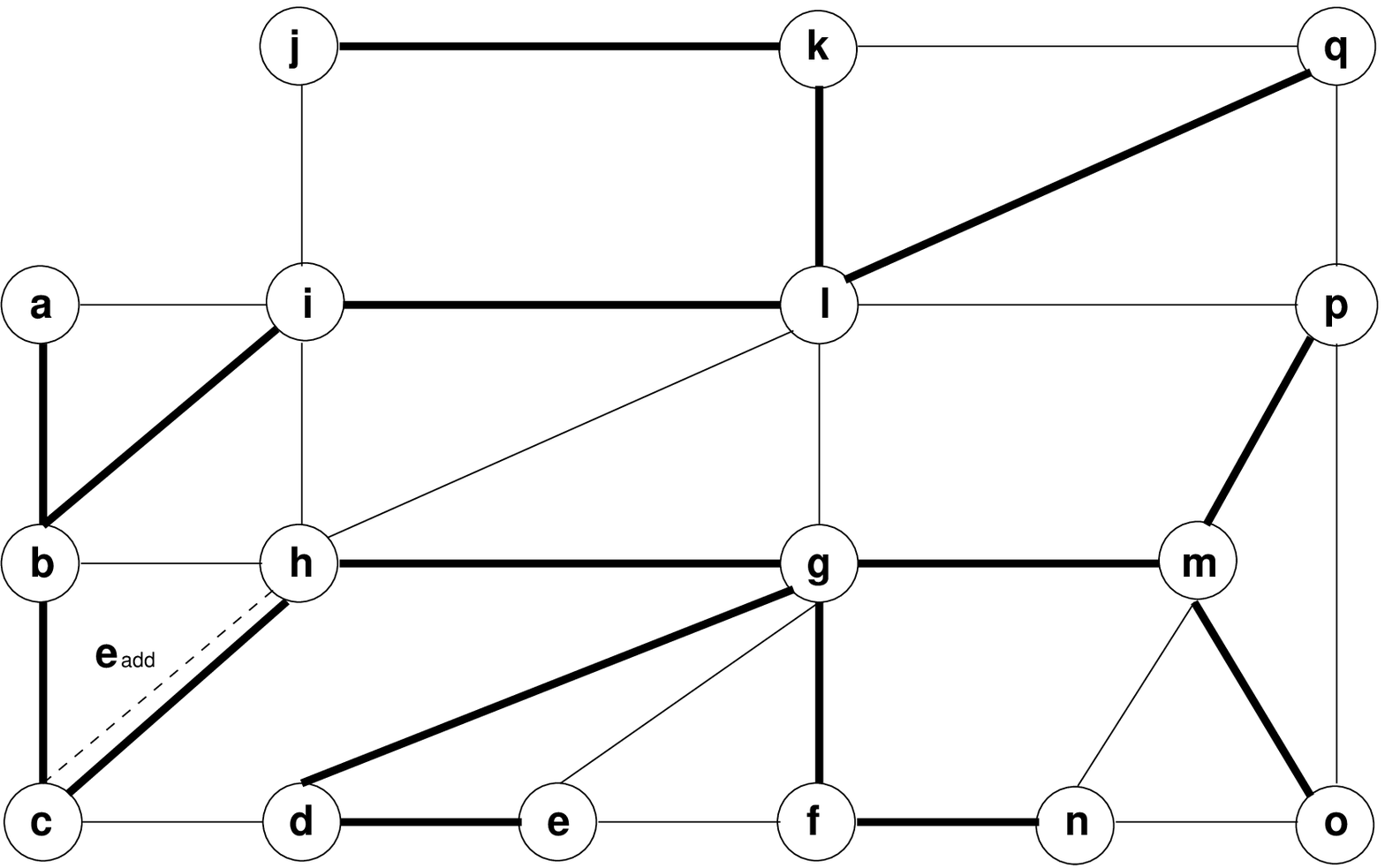}}
\caption{$G^*$ at the end of step $i=2$.}
\end{figure}

\begin{figure}
\centerline{\includegraphics[width=2.5in]{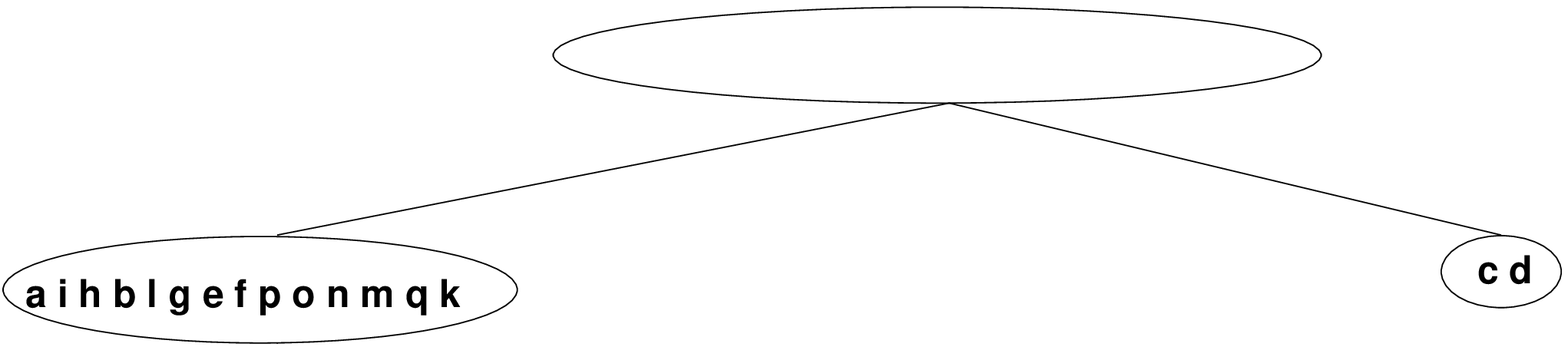}}
\caption{$H$ at the end of step $i=2$.}
\end{figure}

\begin{figure}
\centerline{\includegraphics[width=3.5in]{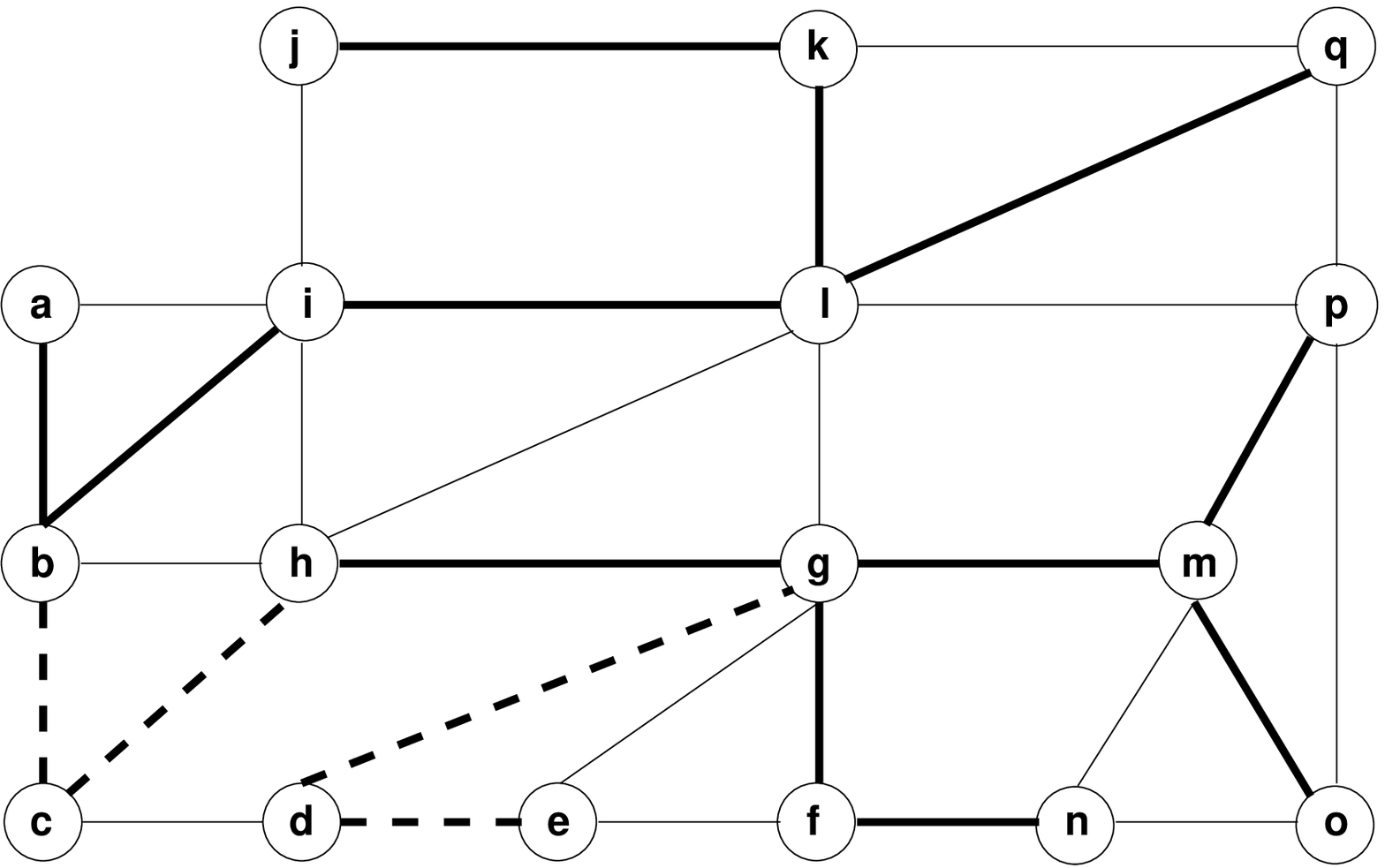}}
\caption{$G^*$ at the end of step $i=3$.}
\end{figure}

\begin{figure}
\centerline{\includegraphics[width=3.5in]{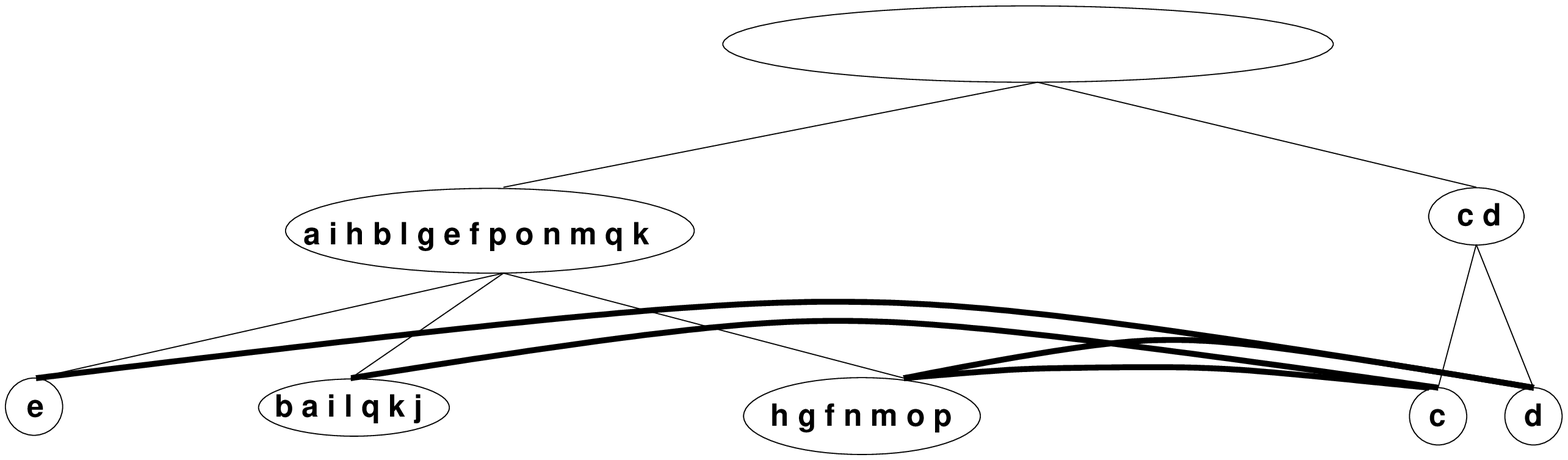}}
\caption{$H$ at the end of step $i=3$.}
\end{figure}

\begin{figure}
\centerline{\includegraphics[width=3.5in]{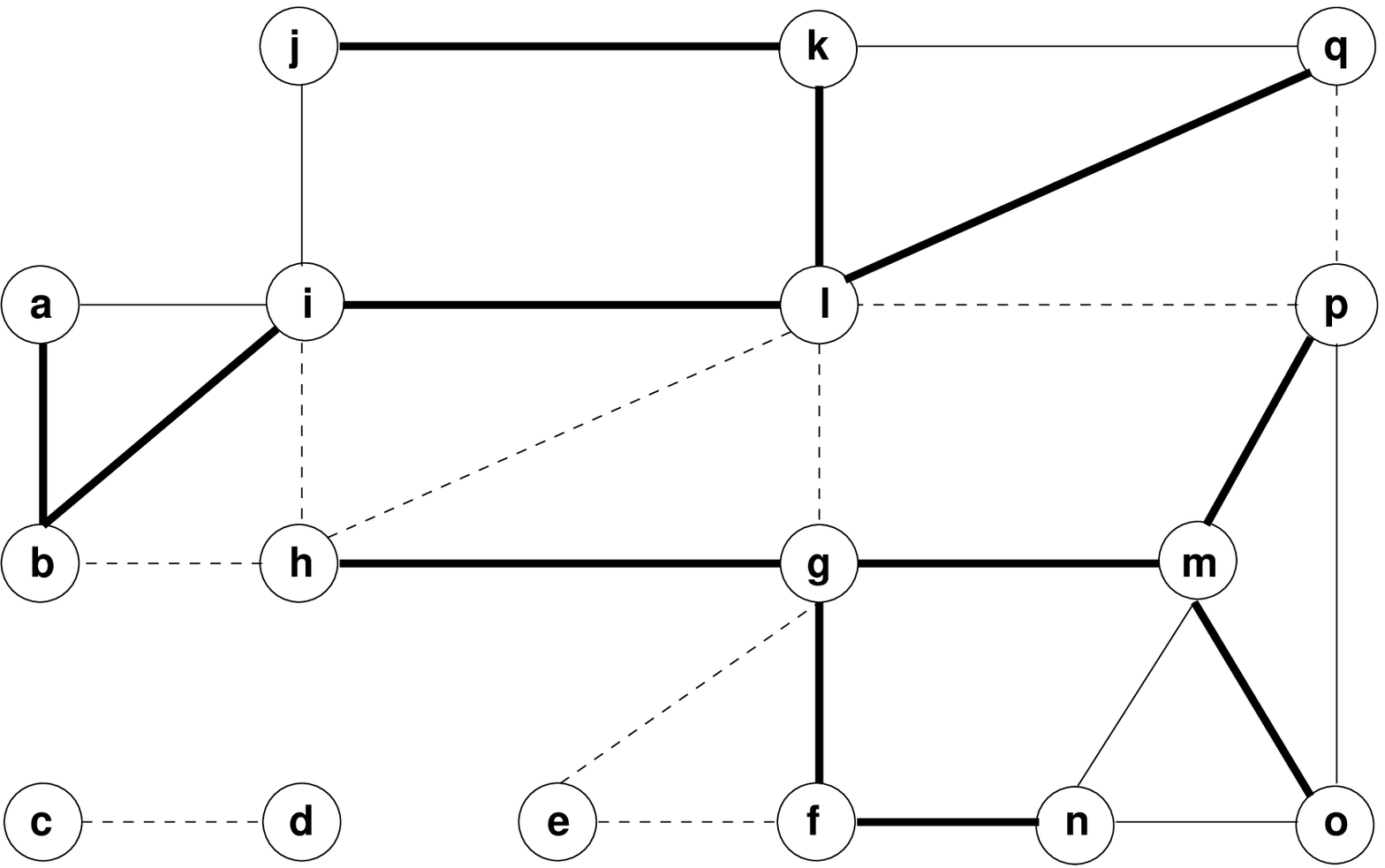}}
\caption{$G^*$ at the end of step $i=4$.}
\label{graph_step3}
\end{figure}


\begin{figure}
\centerline{\includegraphics[width=3.5in]{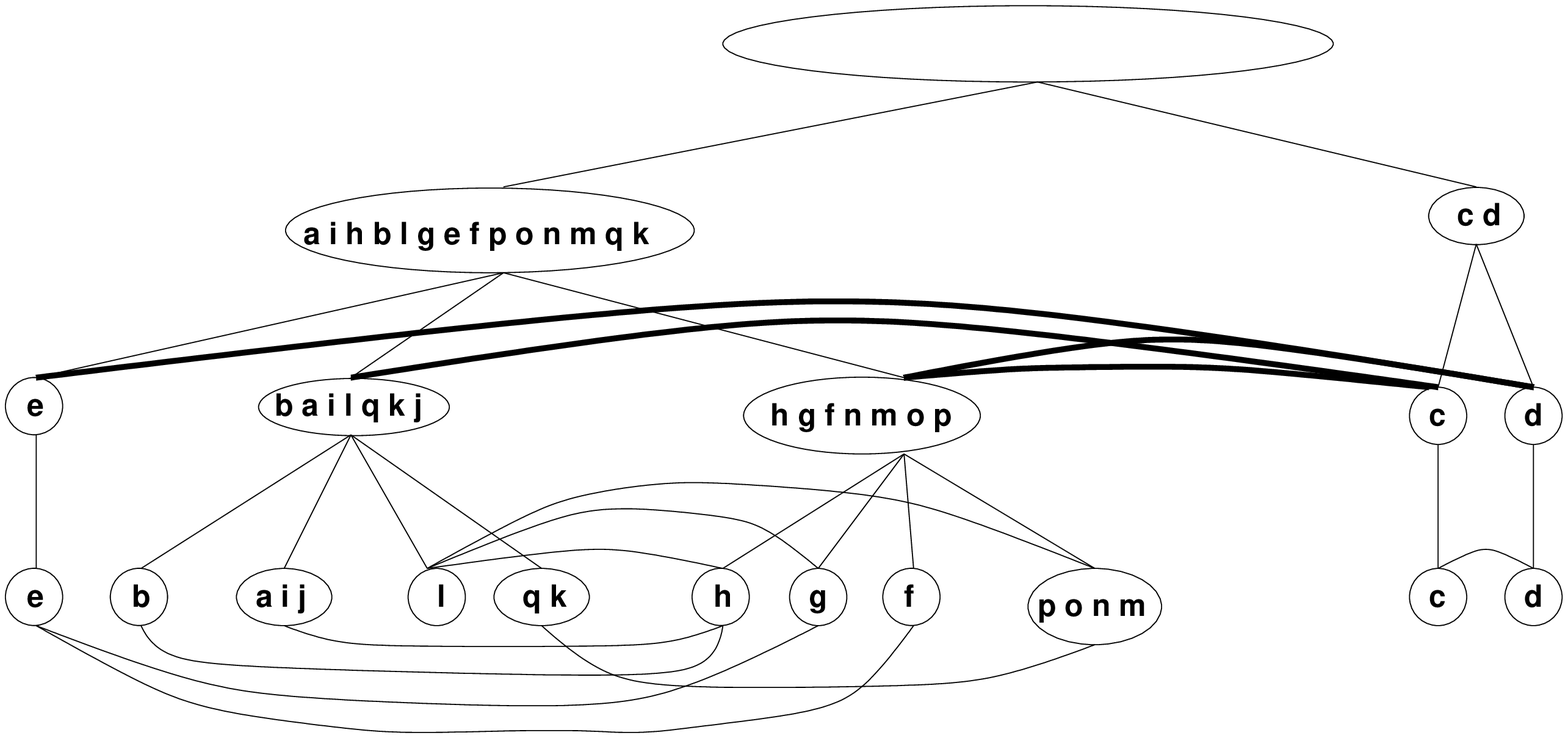}}
\caption{$H$ at the end of step $i=4$.}
\label{hierarchy_step4}
\end{figure}

\begin{figure}
\centerline{\includegraphics[width=3.5in]{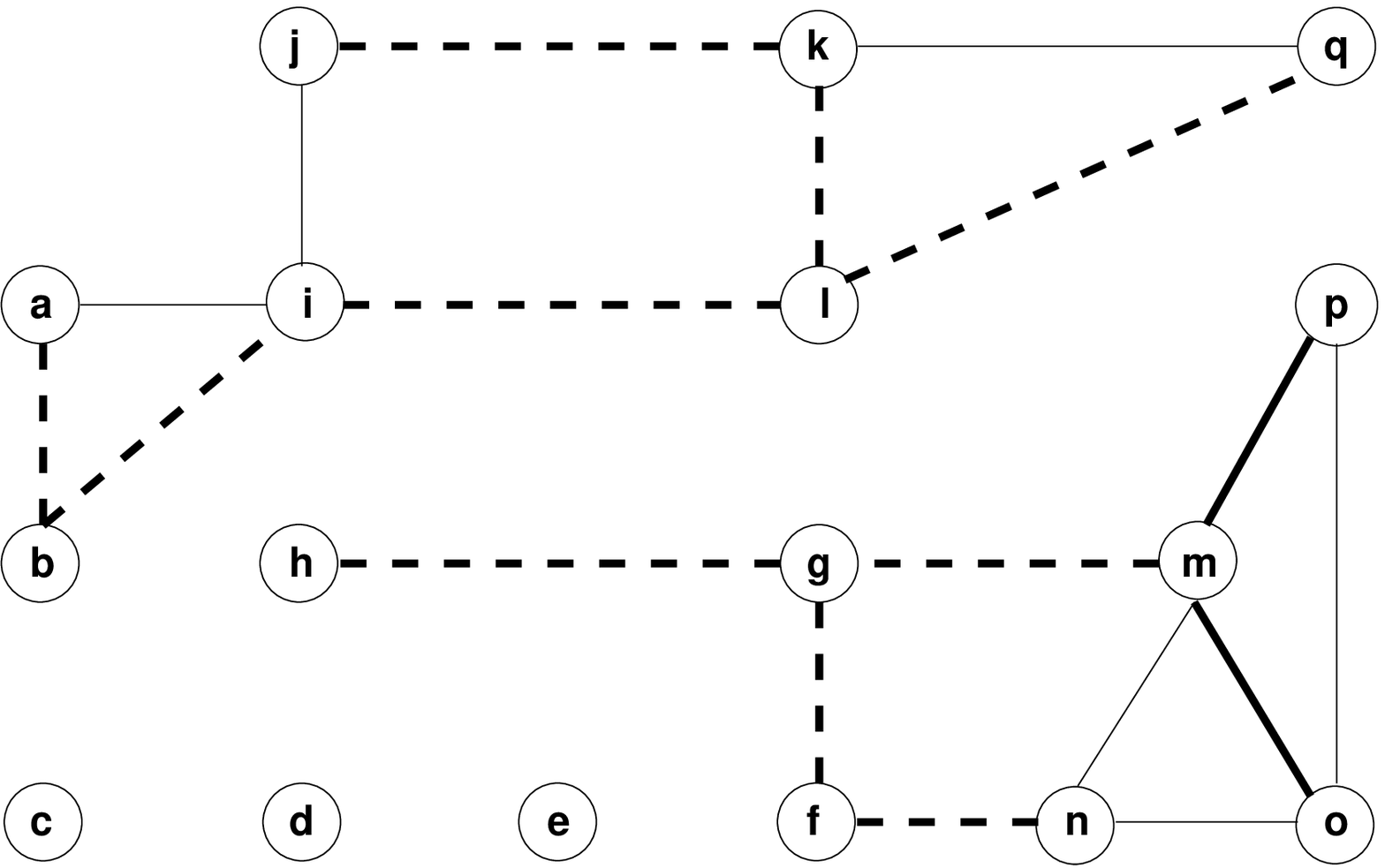}}
\caption{$G^*$ at the end of step $i=5$.}
\end{figure}

\begin{figure}
\centerline{\includegraphics[width=3.5in]{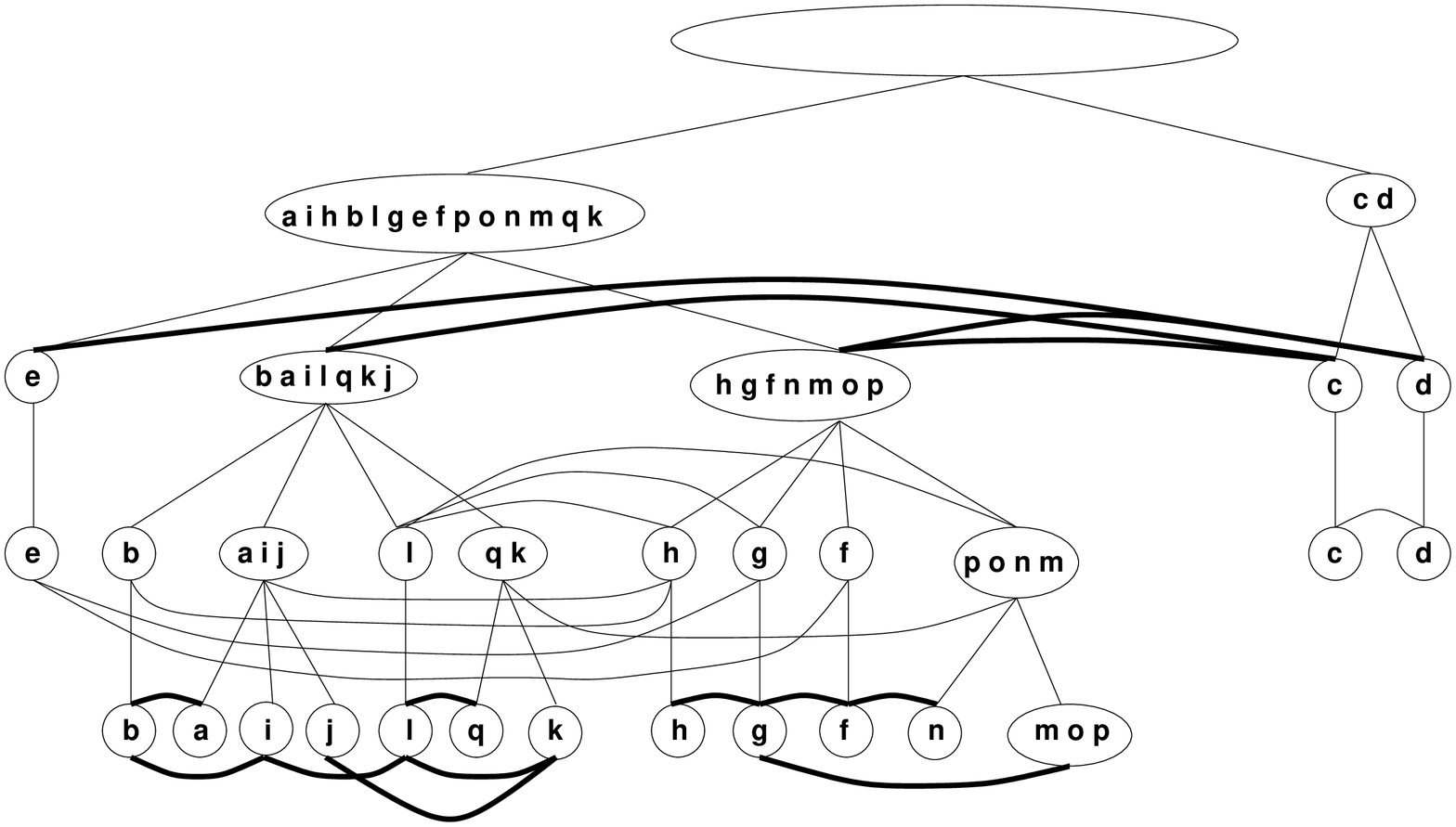}}
\caption{$H$ at the end of step $i=5$.}
\end{figure}

\begin{figure}
\centerline{\includegraphics[width=3.5in]{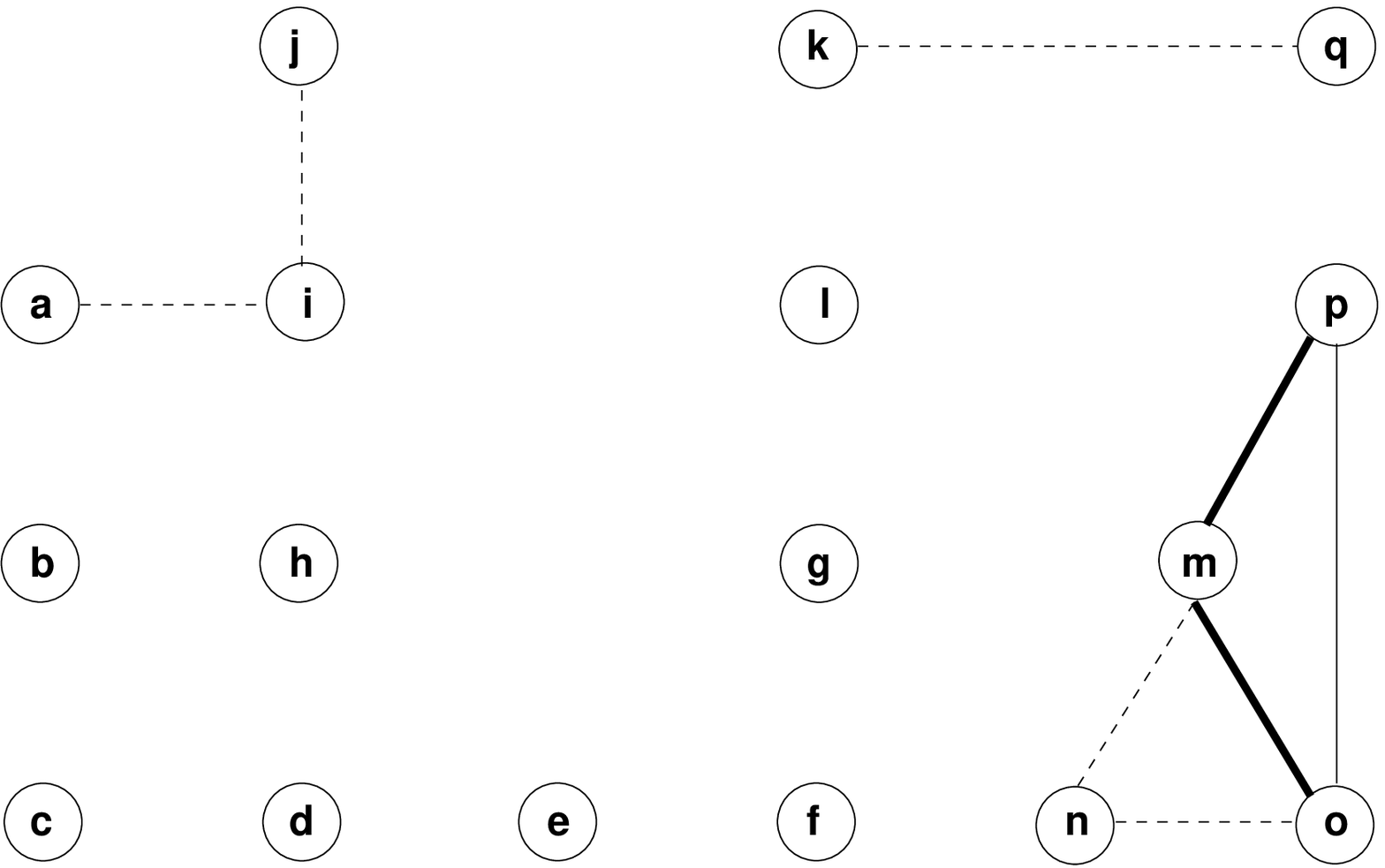}}
\caption{$G^*$ at the end of step $i=6$.}
\end{figure}

\begin{figure}
\centerline{\includegraphics[width=3.5in]{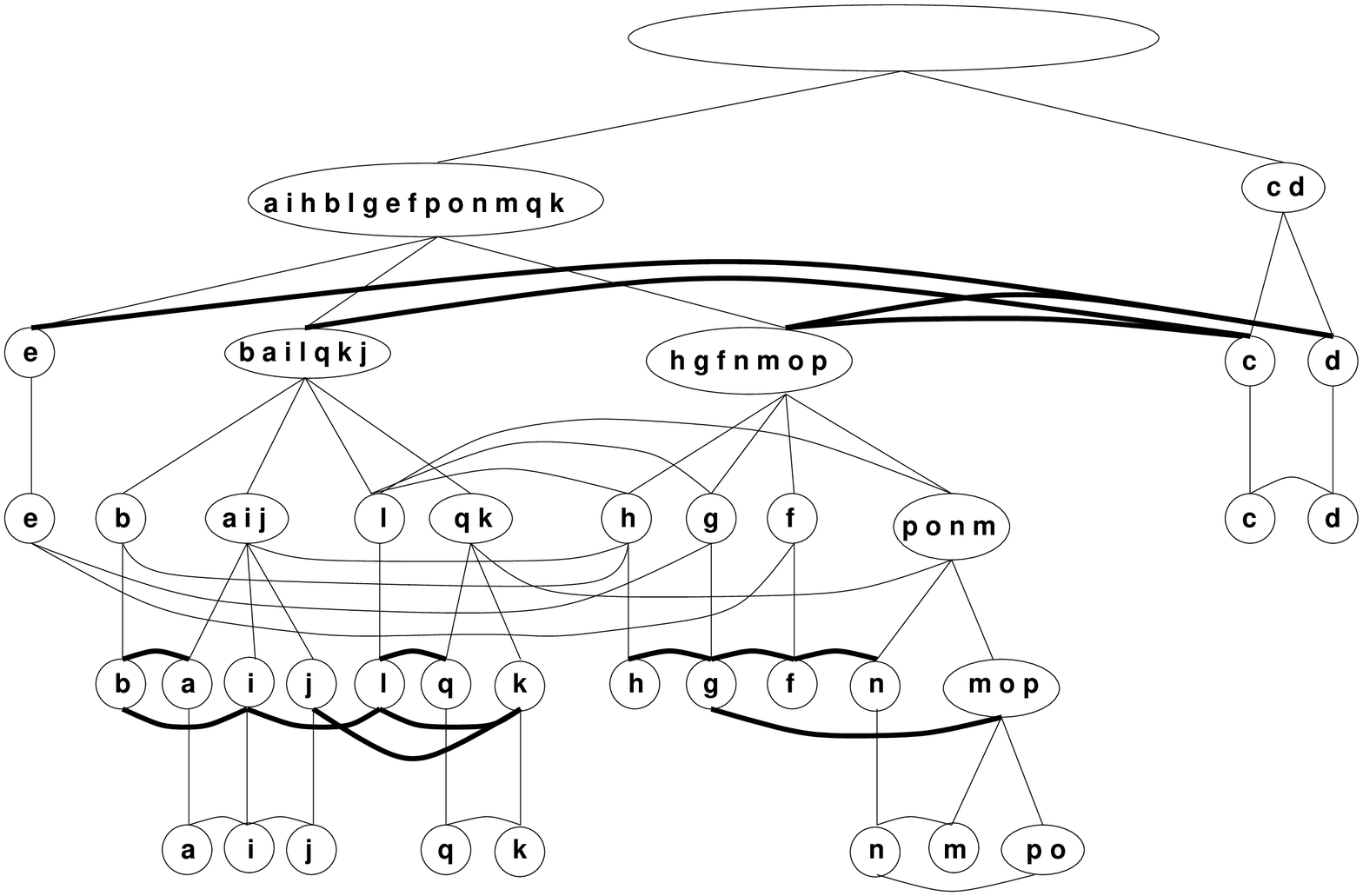}}
\caption{$H$ at the end of step $i=6$.}
\end{figure}
\begin{figure}
\centerline{\includegraphics[width=3.5in]{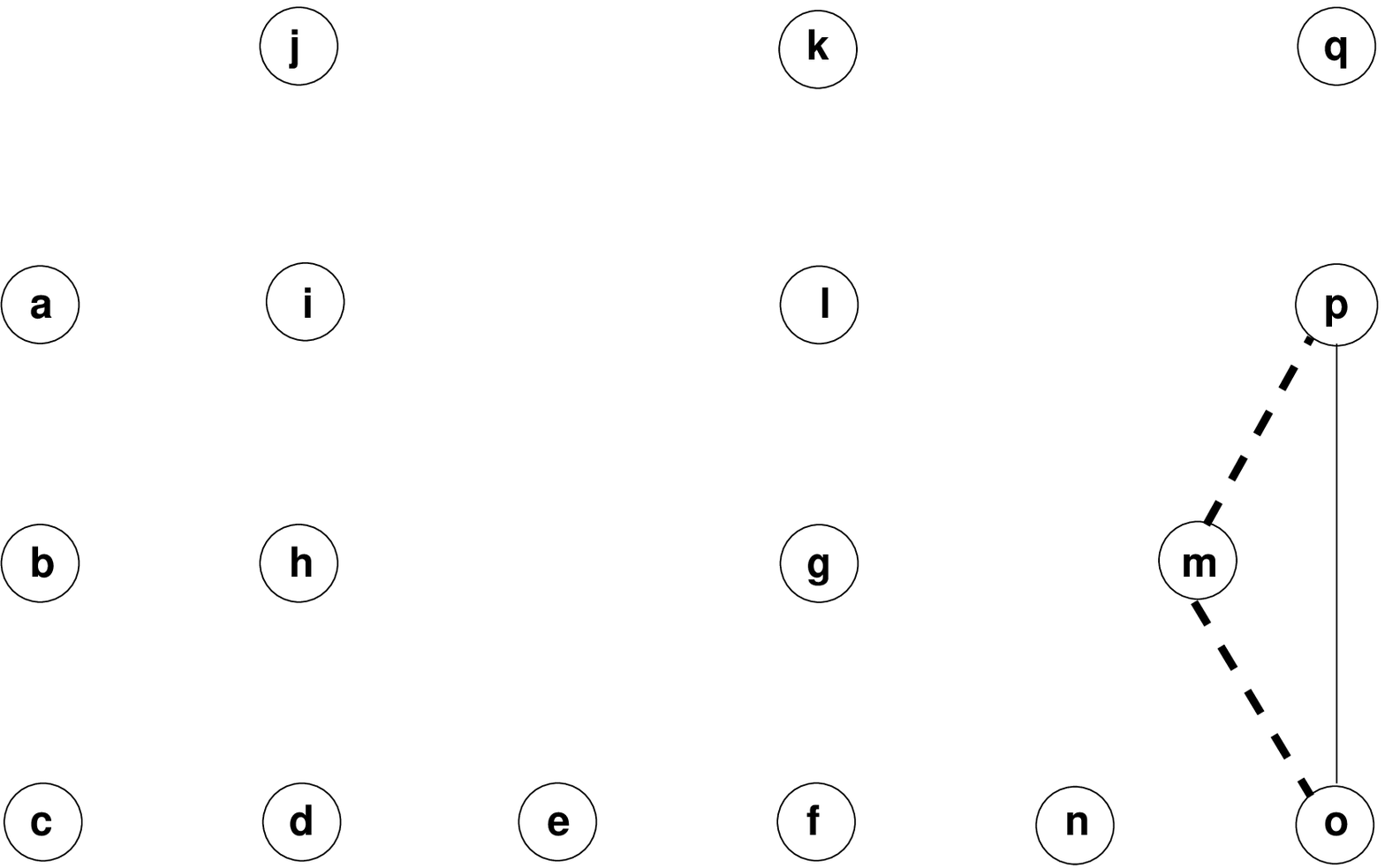}}
\caption{$G^*$ at the end of step $i=7$.}
\end{figure}

\begin{figure}
\centerline{\includegraphics[width=3.5in]{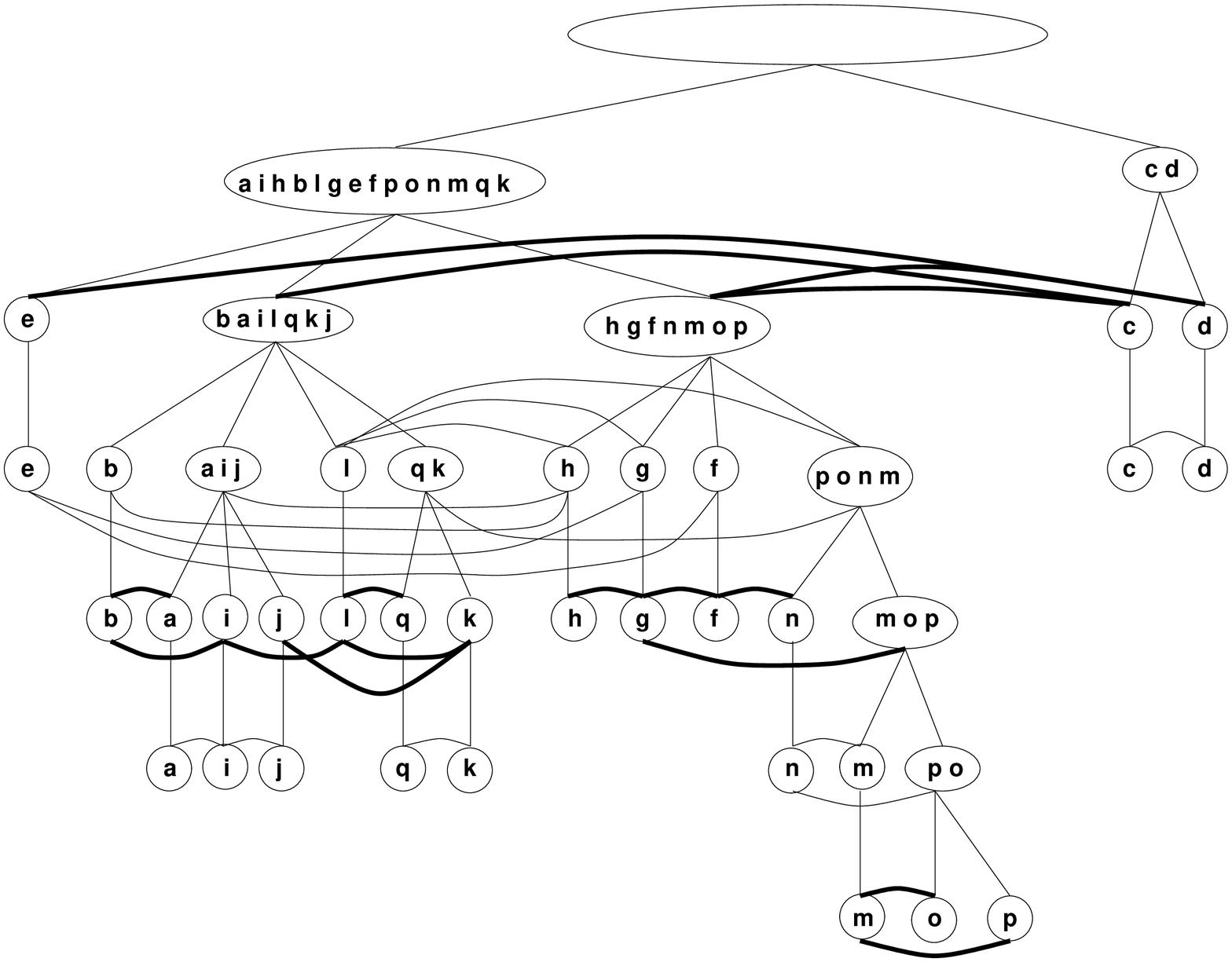}}
\caption{$H$ at the end of step $i=7$.}
\end{figure}

\end{document}